\documentclass[11pt,a4paper]{article}
\usepackage{amsmath,amsfonts,latexsym,amssymb,chicago,annals3,graphicx,here}
\advance \topmargin by -\headheight \advance \topmargin by
-\headsep

\def\m{\mu}

\def\H{{\mathbb H}}
\def\IK{I\!\!K}

\newcommand{\mycite}[1]{{\small \sc \citeNP{#1}}}

\def\a{\alpha}
\def\b{\beta}
\def\R{\mathbb R}

\def\E{{\mathbb E}}
\def\F{{\mathbb F}}

\def\H{{\mathbb H}}

\def\g{\gamma}

\def\labda1{\lambda_1}
\def\labda2{\lambda_2}
\def\m{\mu}

\def\s{\sigma}

\def\comment#1{\relax}

\def\=in{\mathop{\rm =}}

\def\eop{\hfill\mbox{$\Box$}\newline}

\newtheorem{theorem}{Theorem}[section]

\newtheorem{lemma}{Lemma}[section]

\newtheorem{remark}{Remark}[section]

\begin{document}
\title{Isotonic $L_2$-projection test for local monotonicity of a hazard}
\author{Piet Groeneboom and Geurt Jongbloed}
\date{\today}
\affiliation{Delft University of Technology}
\AMSsubject{Primary: 62G10, secondary  62N05.}
\keywords{failure rate, nonparametric testing, bootstrap, isotonic regression}
\maketitle

\begin{abstract} We introduce a new test statistic for testing the null hypothesis that the sampling distribution has an increasing hazard rate on a specified interval $[0,a]$. It is based on a comparison of the empirical distribution function with an isotonic estimate, using the restriction that the hazard is increasing, and measures the excursions of the empirical distribution above the isotonic estimate, due to local non-monotonicity. It is proved in the companion paper \mycite{geurt_piet:11a} that the test statistic is asymptotically normal if the hazard is strictly increasing on the interval $[0,a]$ and certain regularity conditions are satisfied. We discuss a bootstrap method for computing the critical values and compare the test, thus obtained, with other proposals in a simulation study.
\end{abstract}

\section{Introduction}
\label{section:intro}
\setcounter{equation}{0}
In reliability theory and medical statistics, one is often interested in the life time distribution of a certain subject, e.g.\  the distribution of the time it takes before effect of a certain treatment can be noticed or the time it takes before a system device breaks down. In such situations,  it is more natural to model the distribution in terms of its hazard rate (or failure rate) than in terms of the distribution function or density function. Qualitative properties of the hazard rate can be most easily interpreted. These reveal whether or not the device is subject to aging (in case of an increasing failure rate) or not.

Already in the sixties of the preceding century,  estimation of the behavior of the hazard rate  based on a sample from the associated distribution, was studied intensively. The (nonparametric) maximum likelihood estimator (MLE) for the hazard rate is described in, e.g., \mycite{b4:72} and has properties somewhat comparable with the Grenander estimator (MLE) of a decreasing density. Also, procedures were developed to test the null hypothesis of a constant hazard rate (corresponding to an exponential distribution) against the alternative of an increasing hazard (presence of aging). One popular test statistic in this context is the Total Time on Test statistic of  \mycite{prospyke:67}.This is a scale invariant statistic, allowing for efficient computation of Monte Carlo-based critical values.

Only rather  recently, the problem of testing the nonparametric null hypothesis that a hazard rate is (locally) monotone  against the alternative that it is not, has gained  attention. In \mycite{gijbheck:04} local versions of the test statistic of  \mycite{prospyke:67} are studied. In \mycite{durot:08}, the supremum distance between two estimators of the cumulative hazard rate is introduced as test statistic.
In both papers, critical values are obtained using the exponential distribution, which lies `on the boundary of the null hypothesis'. As will be seen in this paper, this choice of the exponential distribution leads to conservative tests when the true underlying distribution has a strictly convex cumulative hazard.

\mycite{hallk:05} developed a test which ``projects" the hazard on the space of nondecreasing hazards by performing a global smoothing of the hazard until it becomes nondecreasing, using kernel estimators. Their criterion for non-convexity of the (non-smoothed) cumulative hazard is compared with the same criterion for a bootstrap sample from the projected (smooth) hazard (this method has been called a ``biased bootstrap"). The idea is that the criterion will be close to zero for bootstrap samples generated from the projected hazard, while the criterion will not be close to zero for the original sample, if the underlying hazard is not monotone. They compare the criterion in the original sample with, say, the 90th or 95th percentile of the distribution of the criterion in the bootstrap samples, and reject the hypothesis of monotonicity if the criterion in the original sample exceeds the chosen percentile of the distribution of the criterion in the bootstrap samples. There are some difficulties with this interesting idea, having to do with non-conservative behavior of this procedure. We will discuss this below.

In this paper we propose another approach, where the type of projection is different from the projection used by \mycite{hallk:05}. As in \mycite{durot:08}, our test statistic is a distance between two estimators for the cumulative hazard. One under the local monotonicity hypothesis and another nonparametric estimator that does not require this monotonicity. Our distance measure is of integral type rather than the supremum distance considered in \mycite{durot:08}. In order to obtain critical values for this test statistic, we propose a bootstrap procedure. Our approach will be described in section \ref{sec:ownapp}. As in the method used by \mycite{hallk:05} and in a certain sense also the method used by \mycite{durot:08}, we will use a bootstrap method for obtaining critical values. In generating the bootstrap samples, we use certain results in \mycite{geurt_piet:11b}, and in the justification of this bootstrap method, we will heavily rely on results in \mycite{geurt_piet:11a}.
Section \ref{sec:simul} contains a simulation study of the various testing procedures, showing that the proposed test has a rather good power, without exhibiting the extreme anti-conservative behavior, exhibited by the method, proposed in \mycite{hallk:05}. The appendix provides the proofs of certain results in section \ref{sec:ownapp}.

\section{Setting and testing procedure}
\label{sec:ownapp}
Consider a sequence of i.i.d.\ random variables $X_1,X_2,\ldots$ with density function $f_0$ on $[0,\infty)$. Denote the distribution function, hazard function and cumulative hazard function associated with $f_0$ by $F_0$, $h_0$ and $H_0$ respectively and recall the relations between these functions:
\begin{equation}
\label{eq:relations}
h_0(x)=\frac{f_0(x)}{1-F_0(x)}, \,\, H_0(x)=-\log(1-F_0(x)),\,\,F_0(x)=1-\exp(-H_0(x)).
\end{equation}
In this paper, we consider the problem of testing local monotonicity (we restrict ourselves to the increasing case; the case of locally decreasing hazard can be considered analogously) of $h_0$. More precisely, given an interval $[a,b]\subset[0,\infty)$, we wish to test
$$
H_{[a,b]}\,:\forall x,y\in[a,b]\,\mbox{ with }x\le y,\, h_0(x)\le h_0(y)
$$
against the alternative that this monotonicity does not hold. Our test statistic is defined as a distance between two estimators for the cumulative hazard function: one under $H_{[a,b]}$ and one that is not.

An estimator for $H_0$ without assuming $H_{[a,b]}$, is just the empirical cumulative hazard function obtained by plugging in the empirical distribution function of the sample $X_1,\dots,X_n$, $\F_n$, in (\ref{eq:relations}):
\begin{equation}
\label{eq:empcumhaz}
\H_n(x)=\left\{
\begin{array}{lll}
-\log\left\{1-\F_n(x)\right\},\,&x\in\left[0,X_{(n)}\right),\\
\infty,\,&x\ge X_{(n)}.
\end{array}
\right.
\end{equation}
Our estimator of $h_0$ under $H_{[a,b]}$ is the least squares estimator, minimizing
the function
\begin{equation}
\label{minim_criterion}
h\mapsto \frac12\int_{[a,b]}h(x)^2\,dx-\int_{[a,b]} h(x)\,d\H_n(x)
\end{equation}
over all nondecreasing functions $h$ on $[a,b]$.

The solution of the problem of minimizing (\ref{minim_criterion}) under the null hypothesis can be constructed explicitly. On $[a,b]$ it is given by the right-continuous derivative of the convex minorant (GCM) of the empirical cumulative hazard function given by (\ref{eq:empcumhaz}), restricted to $[a,b]$. The estimator $\hat H_n$ of $H_0$ under the null hypothesis $H_{[a,b]}$ is therefore defined by
\begin{equation}
\label{eq:Hnhat}
\hat H_n(x)=\left\{\begin{array}{ll}\H_n(x) & x\in[0,a)\cup[b,\infty)\\ {\rm GCM}(y\mapsto\H_n(y)\,:\,a\le y\le  b)(x)& x\in[a,b].\end{array}\right.
\end{equation}
Note that this estimator is continuous at $a$ and, if $b\neq X_i$ for all $i$, at $b$.

Our test statistic for testing the null hypothesis of monotonicity of the hazard on the interval $[a,b]\subset(0,\infty)$ is defined by
\begin{equation}
\label{eq:defTn}
T_n=\int_{[a,b]}\bigl\{\F_n(x-)-\hat F_n(x)\bigr\}\,d\F_n(x).
\end{equation}
where $\hat F_n$ is the distribution function corresponding to $\hat H_n$:
$$
\hat F_n(x)=1-e^{-\hat H_n(x)}.
$$
Note that $T_n\ge0$, since $\hat H_n$ is the greatest convex minorant (hence a minorant) of $\H_n$ on $[a,b]$. Also note that under the alternative hypothesis, $T_n$ will tend to be higher than under the null hypothesis.

\begin{figure}[!ht]
\begin{center}
\includegraphics[scale=0.5]{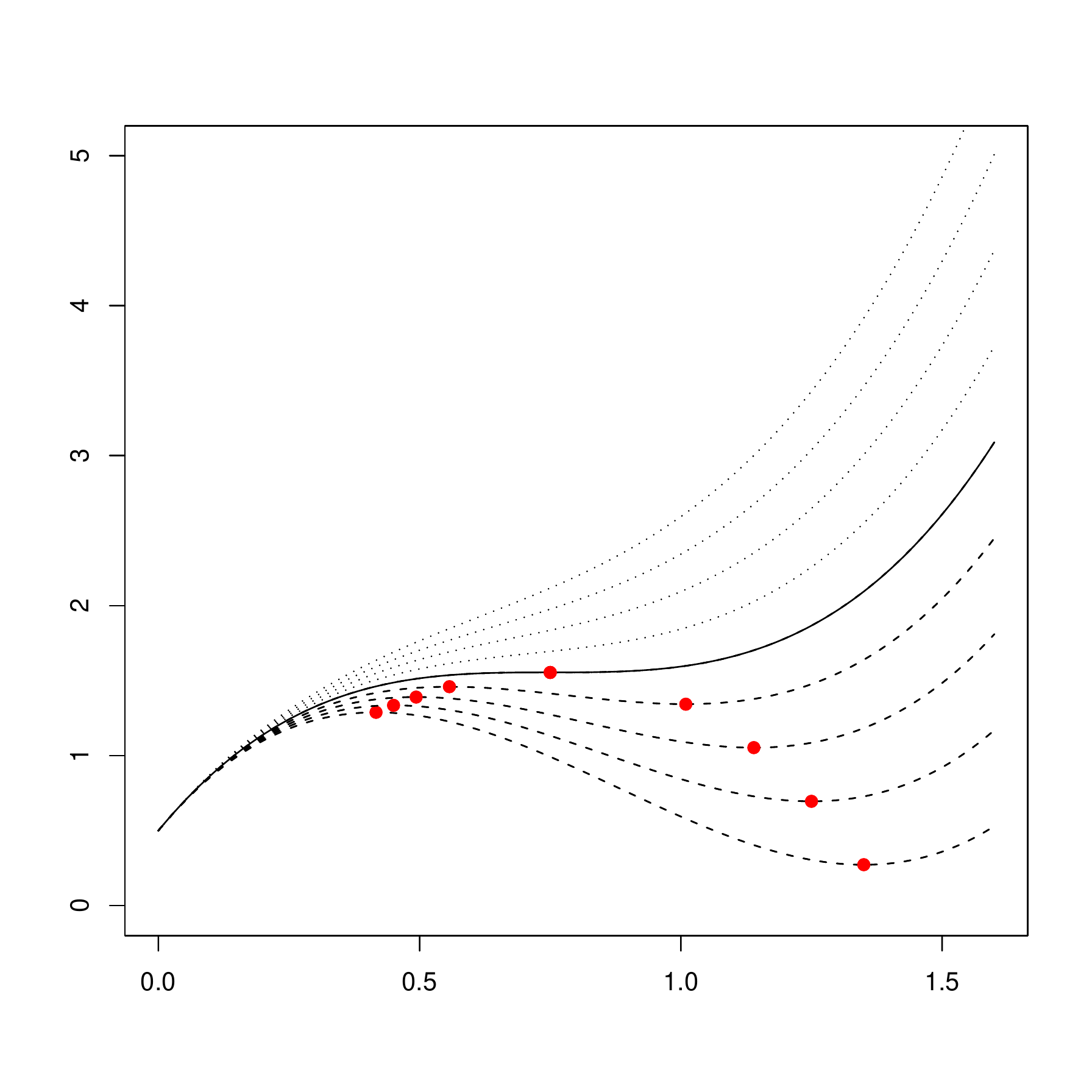}
\end{center}
\caption{The  hazard functions $h^{(d)}$ for $d=-1,-0.75,-0.50,-0.25$ (dashed), $d=0$ (full curve) and $d=0.25,0.50,0.75,1$ (dotted) corresponding to distribution functions (\ref{F_d}). The stationary points are shown by the red dots.}
\label{fig:dfamilyhaz}
\end{figure}

To illustrate the behavior of the estimator we introduce the family of hazards $\{h^{(d)}\,:\,d\in[-1,1]\}$, also considered in \mycite{hallk:05}:
\begin{equation}
\label{h_d}
h_d(x)=\tfrac12+\tfrac52\left\{\left(x-\tfrac34\right)^3+\left(\tfrac34\right)^3\right\}+dx^2,\,x\ge0.
\end{equation}
The corresponding distribution functions are given by:
\begin{equation}
\label{F_d}
F^{(d)}(x)=1-\exp\left\{-\tfrac12x-\tfrac52\left\{\tfrac14\left(x-\tfrac34\right)^4+\left(\tfrac34\right)^3x\right\}-\tfrac13dx^3+\tfrac58\left(\tfrac34\right)^4\right\},\,x\ge0.
\end{equation}
If $d>0$ we get a strictly increasing hazard; if $d<0$, the hazard is decreasing on the interval
$$
\left(\frac34-\frac2{15}d-\frac{2}{15}\sqrt{d^2-\frac{45}{4}d},\frac34-\frac2{15}d+\frac{2}{15}\sqrt{d^2-\frac{45}{4}d}\right)
$$
and if $d=0$ the hazard has a stationary point at $x=3/4$. See Figure \ref{fig:dfamilyhaz} for some hazards in this family.

\begin{remark}
{\rm
Note that we need the constant $\tfrac58\left(\tfrac34\right)^4$ in the exponent to make the distribution function zero at the left endpoint $0$, but that this constant is missing in the formula given below (4.1) on p.\ 1121 in \mycite{hallk:05}.
}
\end{remark}

\begin{figure}[!ht]
\begin{center}
\includegraphics[scale=0.45]{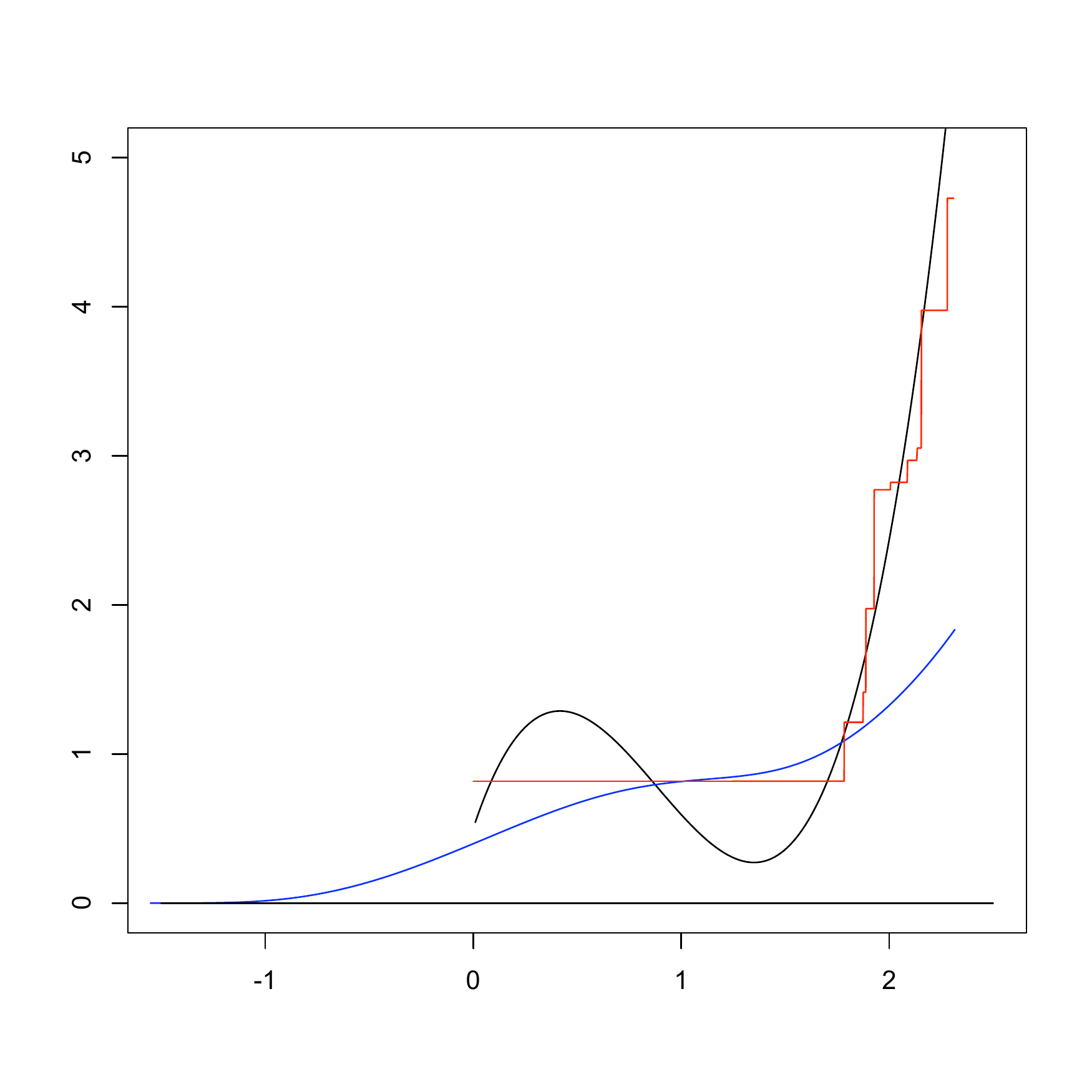}
\includegraphics[scale=0.45]{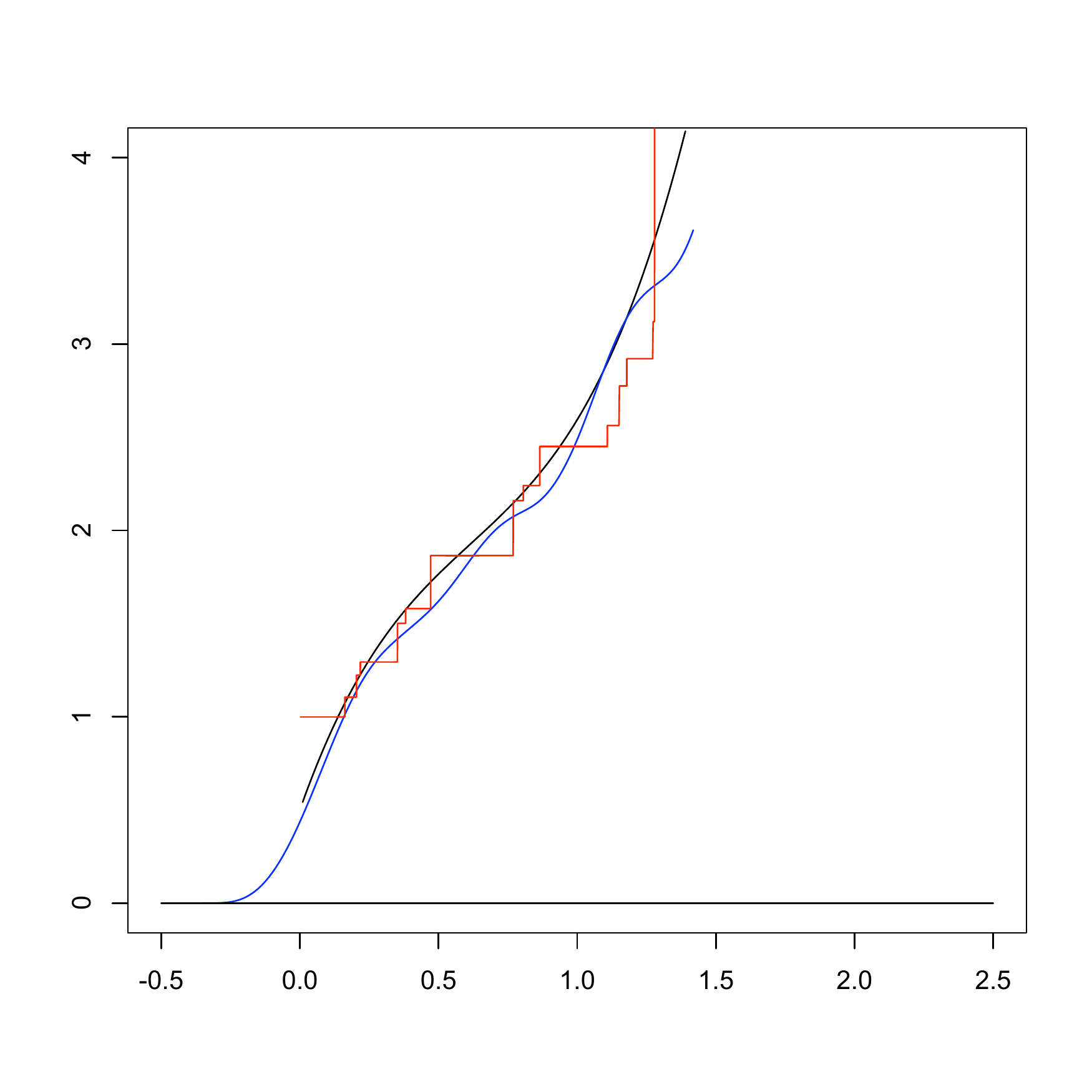}
\end{center}
\caption{The real hazard function $h^{(d)}$ (black), the isotonic estimate $\hat h_n^{(d)}$ of the hazard (red), and the Hall and van Keilegom estimate (blue) of the hazard (after calibration), for a sample of size $n=1000$ from the distribution function $F^{(d)}$. The left panel corresponds to $d=-1$, the right panel with $d=1$.}
\label{fig:3hazards}
\end{figure}

The rather different nature of our isotonic projection of the hazard rate and the projection of \mycite{hallk:05} is illustrated in the left panel of Figure \ref{fig:3hazards}, where $d=-1$.
Their hazard estimate, given by the blue curve in Figure \ref{fig:3hazards} extends (with positive values) to the left of zero and has a slower increase to the right of $2.0$ than the actual hazard which is given by the black curve (which is clearly not monotone). The isotonic projection, on the other hand, only lives on $[0,\infty)$, and follows the steep increase of the real hazard to the right of $2.0$, whereas it only locally corrects for the non-monotonicity. The interval on which the hazard was estimated (and made monotone) was $[0,F^{-1}(0.95))\approx[0,2.31165)$, where $F=F^{(d)}$ in (\ref{F_d}) with $d=-1$.

On the other hand, if we are at the other end of the family $\{F^{(d)}:d\in[-1,1]\}$ at $d=1$, and therefore ``deep inside the null hypothesis region", so to speak, the starting bandwidth for the calibration of the \mycite{hallk:05} method immediately gives an increasing hazard on the interval $[0,F^{-1}(0.95))\approx[0,1.39778)$, where $F=F^{(d)}$ with $d=1$, and the projections of the two methods are less different, see the right panel in Figure \ref{fig:3hazards}.

In order to obtain critical values for statistic $T_n$, there are various possible approaches. The first is to use that its distribution under $H\in H_{[a,b]}$  is stochastically bounded by its distribution under the distribution function with the cumulative hazard function that is obtained by linear interpolation on the interval $(a,b)$
$$
H_{a,b}(x)=\frac{H(b)(x-a)+H(a)(b-x)}{b-a} 1_{[a,b)}(x)+H(x)1_{[0,a)\cup[b,\infty)}(x).
$$
\begin{lemma}
\label{lem:consapp}
For each $H\in H_{[a,b]}$,
\begin{equation}
\label{exponential_ineq}
P_{H}\left(T_n\ge t\right)\le P_{H_{a,b}}\left(T_n\ge t\right)
\end{equation}
for all $t\ge0$.
\end{lemma}

\vspace{0.1cm}

\begin{remark}
{\rm
The proof of Lemma \ref{lem:consapp} (given in the appendix) reveals that the stochastic ordering result also holds if in the definition of $T_n$  $\left(\F_n(x)-\hat F_n(x)\right)^p$ would be used for some $1<p\le\infty$ rather than for $p=1$.
Lemma \ref{lem:consapp} is related to the approximation in  \mycite{durot:08}, section III.  If $H(a)$ and $H(b)$ were known, the distribution of $T_n$ under $H_{a,b}$ could be approximated efficiently using Monte Carlo simulation.
In practice, however, $H(a)$ and $H(b)$ are unknown. In order to really use the approximation, estimates for $H_0$ at $a$ and $b$ are needed. In \mycite{durot:08} this estimation, combined with the stochastic domination of Lemma \ref{lem:consapp}, is called the bootstrap.
}
\end{remark}

It is clear that if the function $H_0$ is {\it strictly} convex on $[a,b]$, the lower bound of Lemma \ref{lem:consapp} may be quite rough. Also, overestimation of the interval $[H_0(a),H_0(b)]$ will lead to a rough bound.  In case of  strict convexity of $H_0$ on $[a,b]$, the convex minorant of its empirical version will tend to wrap tightly around this version whereas in case the cumulative hazard is linear on $[a,b]$ (as in the exponential case), this difference will tend to be bigger. The following theorem, proved in \mycite{geurt_piet:11a}, describes the asymptotic behavior of $T_n$, if the underlying hazard is strictly increasing on $[a,b]$.

\begin{theorem}
\label{th:limitdist}
Let $h_0$ be strictly increasing and positive on the interval $I=[a,b]\subset[0,\infty)$, with a bounded continuous derivative, staying away from zero on $I$. Moreover, let $\zeta(t)$ be the distance at $t$ of the process
$$
W(x)+x^2,\,x\in\R,
$$
to its greatest convex minorant, where $W$ is two-sided Brownian motion, originating from zero.
Then, for $T_n$ as defined in (\ref{eq:defTn}),
$$
n^{5/6}\left\{T_n-ET_n\right\}
\stackrel{{\cal D}}\longrightarrow N\left(0,\s_{F_0}^2\right)
$$
where $N\left(0,\s_{F_0}^2\right)$ is a normal distribution with mean zero and variance $\s_{F_0}^2$, and where
\begin{equation}
\label{mu}
ET_n\sim n^{-2/3}E\zeta(0)\int_a^b \left(\frac{2h_0(t)f_0(t)}{h_0'(t)}\right)^{1/3}\,dF_0(t),
\,n\to\infty,
\end{equation}
and
\begin{equation}
\label{sigma}
\s_{F_0}^2=2\int_0^{\infty}\mbox{\rm covar}(\zeta(0),\zeta(s))\,ds\int_a^{b}\left(\frac{2h_0(t)f_0(t)}{h_0'(t)}\right)^{4/3}\,dF_0(t).
\end{equation}
\end{theorem}

If we want to test the hypothesis that the hazard is strictly increasing on $[a,b]$, we could try to estimate the parameters $\mu$ and $\s$ of Theorem \ref{th:limitdist} and use the limiting normal distribution for the critical values. The difficulty with this approach is that it cannot be used if the derivative of $h_0$ is zero, as is the case, for example, if the underlying distribution is the exponential distribution. For the latter situation we have the following result, proved in \mycite{geurt_piet:11a}.

\begin{theorem}
\label{th:df_theorem2}
Let $U$ be given by
$$
U=\int_0^a\bigl\{1-F_0(x)\bigr\}\left\{W\left(\frac{F_0(x)}{1-F_0(x)}\right)-C(x)\right\}\,dF_0(x),
$$
where $W$ is standard Brownian motion on $[0,\infty)$ and $C$ is the greatest convex minorant of
\begin{equation}
\label{BM_for_constant_haz}
x\mapsto W\left(\frac{F_0(x)}{1-F_0(x)}\right),\,x\in[0,a].
\end{equation}
Suppose that the underlying hazard $h_0$ is constant on $[0,a]$. Then:
\begin{align*}
n^{1/2}T_n\stackrel{{\cal D}}\longrightarrow U,\,n\to\infty.
\end{align*}
\end{theorem}

So we see that in this situation the rate of convergence drops from $n^{5/6}$ to $n^{1/2}$, and also that the limit distribution is no longer normal. In the case of the family $h^{(d)}$, we therefore enter a completely different regime of asymptotic behavior when the parameter $d$ passes zero. The most natural way to deal with this difficulty seems to us to be a bootstrap procedure which we will now describe. We will only prove that our bootstrap method works under the conditions of Theorem \ref{th:limitdist}. We conjecture, however, that our method will also work under the conditions of Theorem \ref{th:df_theorem2} (possibly with a slightly modified version of $\tilde H_n$), but this is still an open question.

The proposed method runs as follows. First estimate the cumulative hazard function under the null hypothesis by a smooth estimator, having the property that the corresponding hazard satisfies the null hypothesis. Then draw samples of size $n$ from this estimate $B$ times and compute $B$ times the bootstrap version of the test statistic: $T^*_{n,i}$, $1\le i\le B$. Finally, approximate the distribution of $T_n$ under the true cumulative hazard function $H_0$ (assumed to belong to $H_{[a,b]}$)  by the empirical distribution of these bootstrap values and its critical value at (for example) level $10\%$ by the $90$-th percentile of this generated set of bootstrap values. In fact, \mycite{hallk:05} also use a bootstrap procedure of this type, but based on a totally different ``projection" of the hazard estimate on the set of increasing hazards.

In the further description of the method, we take the left endpoint of the interval at the origin and denote the right endpoint by $a$, as in \mycite{geurt_piet:11a}.
In order to prevent inconsistency of $\hat h_n$ at the endpoints, we define a penalized version of $\hat h_n$, $\hat{h}_n^{[p]}$, as the derivative of the penalized cusum diagram consisting of the points
\begin{equation}
\label{cusum2}
(0,0),\qquad \left(X_{(i)},\H_n(X_{(i)}-)+2n^{-2/3}\right),\,X_{(i)}<a,\qquad \left(a,\H_n(a-)\right).
\end{equation}
The left derivative of the present cusum diagram minimizes the criterion
\begin{equation}
\label{LS_criterion2}
\tfrac12\int_0^a h(x)^2\,dx-\int_{[0,a]} h(x)\,d\H_n(x)-\a_n h(0)+\b_n h(a),
\end{equation}
where $\a_n=\b_n=2n^{-2/3}$, over all nondecreasing functions $h$ on $[0,a]$. The reason for choosing a penalty of order $cn^{-2/3}$ is explained in \mycite{geurt_piet:11b}, where also a proof of the consistency of
this estimator at the boundary points is given.

For $x\in[0,a]$, we estimate the hazard by kernel smoothing of $\hat h_n^{[p]}$. Let $K$ be the triweight kernel
\begin{equation}
\label{kernel}
K(u)=\tfrac{35}{32}\left\{1-u^2\right\}^31_{[-1,1]}(u),\,u\in\R.
\end{equation}
This is a mean zero probability density with second moment $1/9$.
Then, define for bandwidth $b_n>0$
\begin{equation}
\label{smooth_haz_est}
\tilde h_n(x)=\int K_{b_n}(x-y)\,d\hat H_n^{[p]}(y)=\int K_{b_n}(x-y)\,\hat h_n^{[p]}(y)\,dy,
\end{equation}
where $K_{b_n}(u)=K(u/b_n)/b_n$.
Equation (\ref{smooth_haz_est}) can then be written as
\begin{align*}
\tilde h_n(x)&=\int K_{b_n}(x-y)\,\int_{0}^y d\hat h_n^{[p]}(u)\,dy
=\iint_{u<y} K_{b_n}(x-y)\,dy\,d\hat h_n^{[p]}(u)\\
&=\int_{u=0}^{x+b_n}\int_{y=u}^{x+b_n}K_{b_n}(x-y)\,dy\,d\hat h_n^{[p]}(u)
=\int_{u=0}^{x+b_n}\IK\left(\frac{x-u}{b_n}\right)\,d\hat h_n^{[p]}(u),
\end{align*}
where
$$
\IK(u)=\int_{-\infty}^uK(w)\,dw=1_{[-1,1)}(u)\int_{-1}^u K(w)\,dw+1_{[1,\infty)}(u).
$$
The corresponding estimate of the $h_0^{\prime}$ and $H_0$ are then given by
\begin{equation}
\label{tilde_H_n}
\tilde h_n'(x)=\int K_{b_n}(x-y)\,d\hat h_n^{[p]}(y) \mbox{ and }\tilde H_n(x)=\int_0^x \tilde h_n(u)\,du.
\end{equation}

\begin{figure}[!ht]
\begin{center}
\includegraphics[scale=0.5]{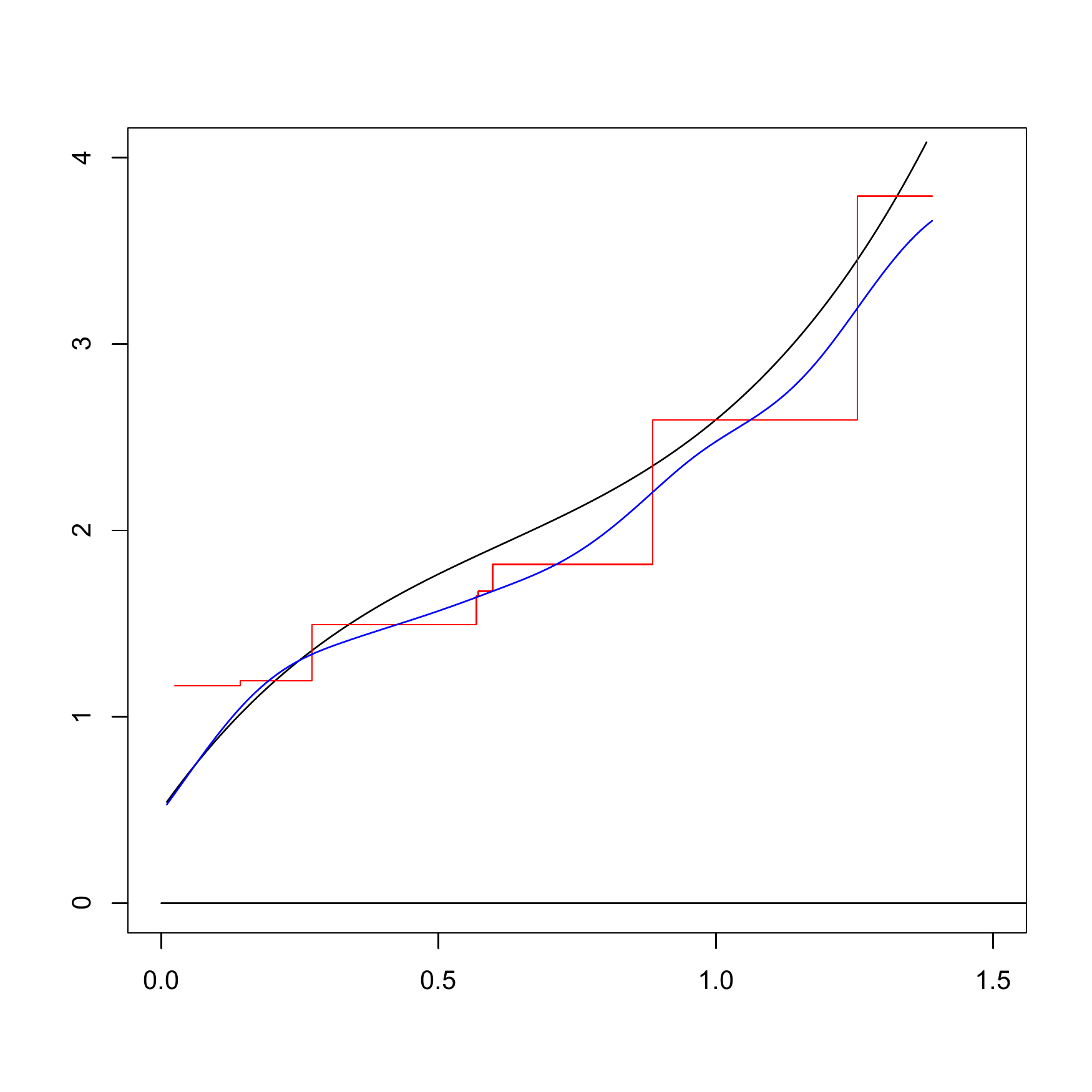}
\end{center}
\caption{The estimate $\tilde h_n$ (blue) of the hazard $h^{(1)}$ (black) of the family $\{h^{(d)}:d\in[-1,1]\}$ for a sample of size $n=100$, together with the (penalized with $2n^{-2/3}\approx0.093$) isotonic estimate $\hat h_n$ (red) on the $95\%$ percentile interval $[0,\left(F^{(1)}\right)^{-1}(0.95)]$. Bandwidth $b_n= n^{-1/4}\approx0.316$.}
\label{fig:hazards1000}
\end{figure}

In justifying this method for testing that the hazard is strictly increasing on $[0,a]$, we use the following bootstrap version of Theorems \ref{th:limitdist}, which will be proved in section \ref{sec:appendix}.

\begin{theorem}
\label{th:bootstrap1}
Let the conditions of Theorem \ref{th:limitdist} be satisfied, and let $\tilde H_n$ be the estimate of the cumulative hazard function under the null hypothesis, defined by (\ref{tilde_H_n}), and based on a sample $X_1,\dots,X_n$ from $F_0$, where we take a vanishing bandwidth $b_n$, satisfying  $b_n\gtrsim n^{-1/4}$. Let $X_1^*,\dots,X_n^*$ be a bootstrap sample generated by $\tilde H_n$ and let $\F_n^*$ and $\hat F_n^*$ be the (bootstrap) empirical distribution function and corresponding estimate $\hat F_n^*$, based on the greatest convex minorant of the function $x\mapsto -\log(1-\F_n^*(x-))$, respectively. Finally, let $T_n^*$ be defined by
$$
T_n^*=\int_{[0,a]}\bigl\{\F_n^*(x-)-\hat F_n^*(x)\bigr\}\,d\F_n^*(x),
$$
and let its (bootstrap) expectation be defined by
$$
E^* T_n^*=\int_{[0,a]}\E\left\{\F_n^*(x-)-\hat F_n^*(x)\bigm|\tilde F_n\right\}\,d\tilde F_n(x),
$$
where $\tilde F_n(x)=1-\exp\{-\tilde H_n(x)\}$. Then we have, almost surely,
$$
n^{5/6}\left\{T_n^*-E^*T_n^*\bigm|X_1,\dots,X_n\right\}\stackrel{{\cal D}}\longrightarrow N\left(0,\s_{F_0}^2\right),
$$
as $n\to\infty$, where  $\s_{F_0}^2$ is given in Theorem \ref{th:limitdist}.
\end{theorem}

\section{A simulation study}
\label{sec:simul}
In this section we compare the power behavior of the test based on our test statistic $T_n$, defined by (\ref{eq:defTn}), with other test statistics for the families that were also considered in \mycite{hallk:05}. The bootstrap resampling for $T_n$ was  done by taking $B=2000$ samples from the estimate $\tilde H_n$ defined in (\ref{tilde_H_n}) with bandwidth $b_n=n^{-1/4}$. For the estimator  $\hat h_n^{[p]}$ on which $\tilde H_n$ is based (see (\ref{cusum2})), the penalty was taken equal to $2n^{-2/3}$ . The sample was generated by first generating a standard exponential sample $E_1,\dots,E_n$, producing the bootstrap sample via
$$
X_i^*=\tilde H_n^{-1}(E_i),\, 1\le i\le n.
$$
In this way, $B$ values $T_n^*$ were obtained.
The critical value is taken to be the $90$th percentile of these values of $T_n^*$.

Below we also make a comparison with a test, proposed in \mycite{durot:08}, referred to as Durot test in the sequel. This test is based on the supremum distance between the empirical cumulative hazard function and its greatest convex minorant:
$$
T_{n,Durot}=\sup_{x\in[0,a]}\{\H_n(x)-\hat H_n(x)\}.
$$
For determining a critical value, again $B=2000$ random standard exponential samples were generated, and the value
$$
T_{n,Durot}^*=\sup_{x\in[0,\H_n(a)]}\{\H_n^*(x)-\hat H_n^*(x)\}
$$
was determined for each such ``bootstrap" sample (taking the interval $[0,\H_n(a)]$ as interval of convexity). The critical value was then taken to be the $90$th percentile of the so obtained values of $T_{n,Durot}^*$. Note that this procedure is equivalent to the procedure that first estimates the (constant) hazard rate on $[0,a]$ by  $\H_n(a)/a$, then takes bootstrap samples from the exponential distribution with this hazard rate and finally determines the supremum distance between the two resulting estimators on the interval $[0,a]$. 

In Table \ref{table:powers_d1}, four tests are compared: the test, based on $T_n$, the test proposed in \mycite{hallk:05} (in the sequel referred to as HvK test), the Durot test and the integral statistic version of this statistic, where we replace the maximum distance statistic by $T_n$, defined by (\ref{eq:defTn}) using Durot's method of approximating the critical value. In this table the tests are compared on the fixed interval $[0,a]=[0,F_0^{-1}(0.95)]$ (instead of on the random interval $[0,\F_n^{-1}(0.95)]$, as in \mycite{hallk:05}. In all cases we generated 2000 samples, and also $B=2000$ bootstrap samples from each original sample.

The simulations for \mycite{hallk:05} took rather long, since repeated density estimation is needed at each step in view of the needed calibration of the bandwidth to create a non-decreasing hazard in the original sample. Also, one has to compute an estimator of the distribution function, the density, and the derivative of the density to check whether one gets a nondecreasing hazard on the chosen interval at the critical bandwidth. The estimation of the density and its derivative was speeded up by using Fast Fourier Transform, and the distribtution function was computed by numerically integrating the density estimate.

\begin{figure}[!ht]
\begin{center}
\includegraphics[scale=0.45]{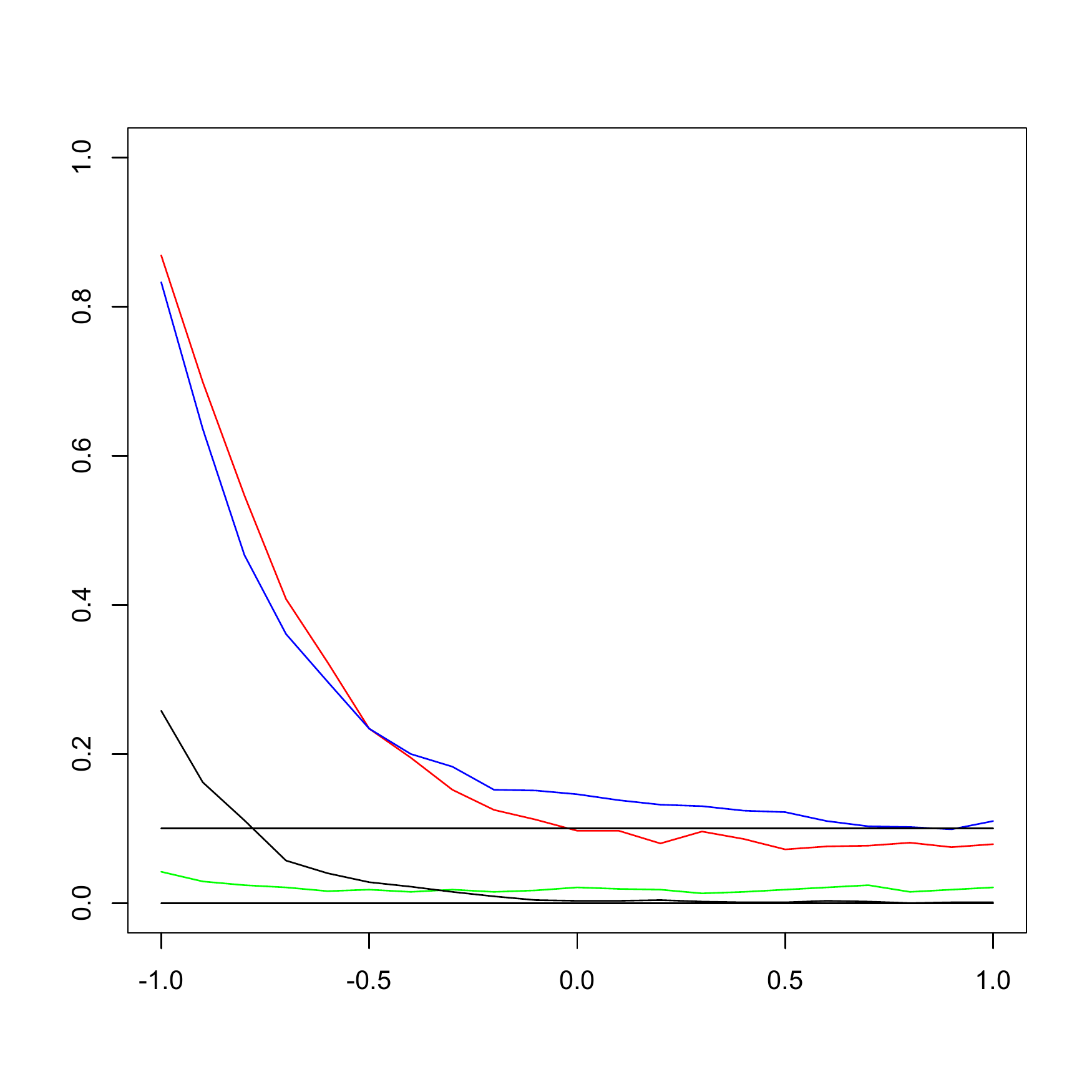}
\includegraphics[scale=0.45]{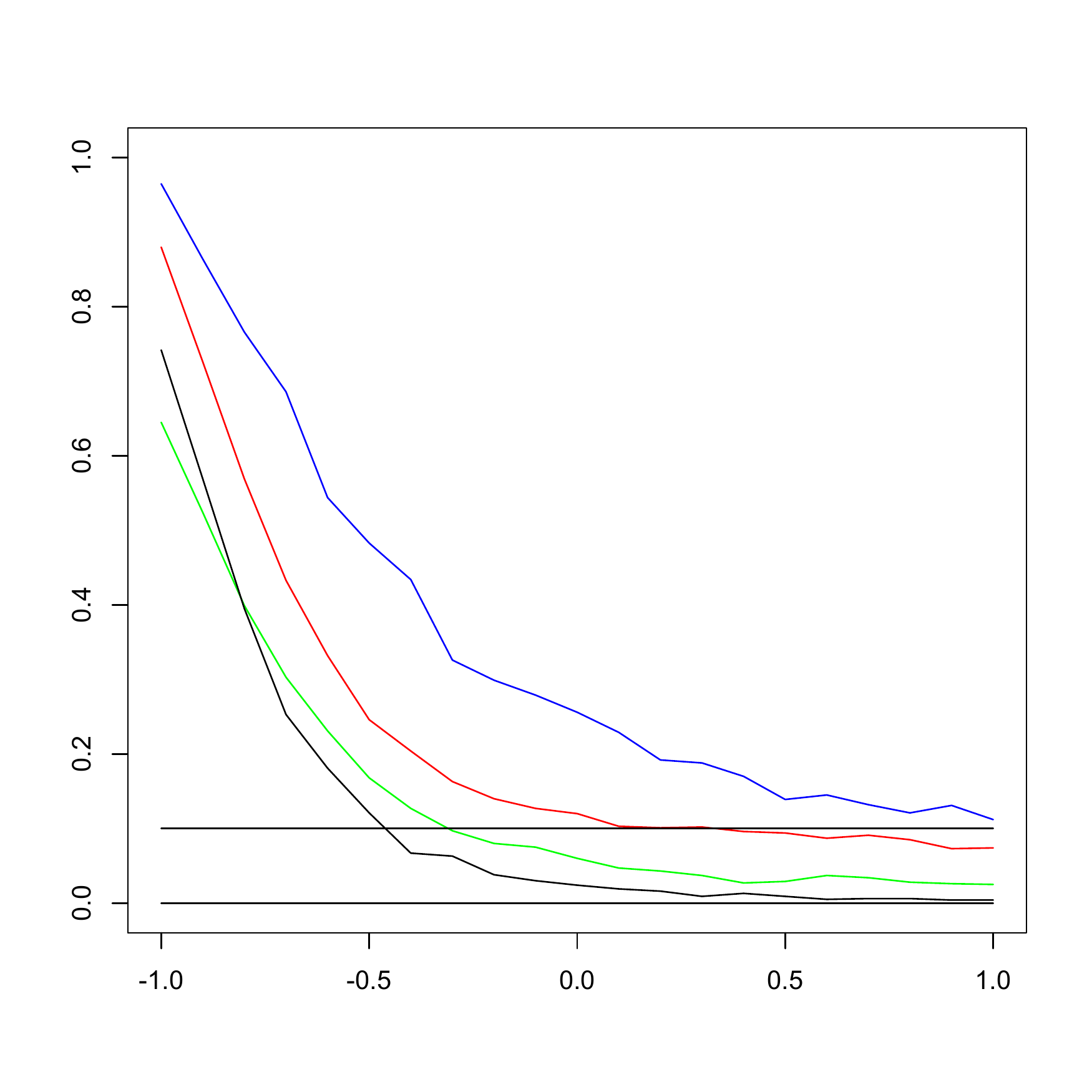}
\end{center}
\caption{The estimated power curves for the family $\{F^{(d)}:d\in[-1,1]\}$, the isotonic test statistic $T_n$, defined by (\ref{eq:defTn}) (red), for the HvK test (blue), the Durot test  (green), and the integral statistic version of this method (black). The sample size $n=50$ and the estimation interval is $[0,F_0^{-1}(0.95)]$ in the left panel, $[0,F_0^{-1}(0.8)]$ in the right panel.}
\label{fig:powers5080_fixed}
\end{figure}

It is seen from Table \ref{table:powers_d1} that the test based on $T_n$ is slightly more powerful for the alternatives $F^{(d)}$ for $d\in[-1,-0.5]$ than the HvK test.
Table \ref{table:powers_d2} shows that the HvK test  is rather anti-conservative. This seems to suggest that the high power in the region $d\in[-0.5,0]$ is at least partly due to the anti-conservative behavior of this test. The Durot test  is very conservative for this interval, as is to be expected, since the estimated critical value is based on the exponential distribution. The test based on $T_n$, has a middle position: it is more conservative than the HvK test but less conservative than the Durot  test. A graphical comparison of the power functions is given in the left panel of Figure \ref{fig:powers5080_fixed}.

Interestingly, the power of the Durot test  increases considerably if we take a smaller interval $[0,F_0^{-1}(0.80))$. In fact, the Durot test proposed  is derived under the assumption that not all order statistics belong to the interval $[0,a]$. But this often happens if we take $[0,a]=[0,F_0^{-1}(0.95))$, in particular for the ``bootstrap samples".

Another reason for the higher power of the Durot test  in this situation is the fact that the isotonic projection of the hazard $h^{(d)}$, for $d\in[-1,0]$ is almost constant on the interval $[0,F_0^{-1}(0.8))$, since we miss the steeply increasing part of the hazard from $F_0^{-1}(0.8)$ to $F_0^{-1}(0.95)$, so sampling from the isotonic projection is almost the same as sampling (locally) from an exponential distribution in this case.

The results are shown in Tables \ref{table:powers_d3} and \ref{table:powers_d4}, and the right panel of Figure \ref{fig:powers5080_fixed}. If one chooses this interval, the HvK test is very powerful, but also {\it very} anti-conservative. For example, for $d=0$ (which belongs to the null hypothesis region) one gets an estimated rejection probability of more than $25\%$ instead of the desired $10\%$.

\begin{table}[h]
\centering
\caption{Estimated powers for model (\ref{F_d}), where $\a=0.1$, $n=50$, and $d=-1,-0.9,\dots,-0.1$. The estimation interval is $[0,F_0^{-1}(0.95)]$.}
\label{table:powers_d1}
\vspace{0.5cm}
\begin{tabular}{|l|c|c|c|c|c|c|c|c|c|c|c|}
\hline
 & \multicolumn{10}{c|}{$d$}\\
\hline
 & $-1$ &  $-0.9$ & $-0.8$ & $-0.7$ & $-0.6$ & $-0.5$ & $-0.4$ & $-0.3$ & $-0.2$ & $-0.1$\\
\hline
$T_n$ &0.869 & 0.699 & 0.547 & 0.408 &0.323 &0.234 &0.195 &0.152 &0.125 &0.112\\
\mbox{HvK} &0.833 & 0.636 & 0.467 & 0.361 &0.297 &0.234 &0.200 &0.183 &0.152& 0.151\\
\mbox{Durot} & 0.042 & 0.029 & 0.024 & 0.021 & 0.016 & 0.018 &0.015 &0.018 &0.015 &0.017\\
\mbox{Durot, $T_n$} & 0.258 & 0.162 & 0.111 & 0.057 & 0.040 & 0.028 &0.022 &0.015 &0.009 &0.004\\
\hline
\end{tabular}
\end{table}

\begin{table}[h!]
\centering
\caption{Estimated rejection probabilities for model (\ref{F_d}) under the null hypothesis, where $\a=0.1$, $n=50$, and $d=0,0.1,\dots,1$. The estimation interval is $[0,F_0^{-1}(0.95))$.}
\label{table:powers_d2}
\vspace{0.5cm}
\begin{tabular}{|l|c|c|c|c|c|c|c|c|c|c|c|c|}
\hline
 & \multicolumn{11}{c|}{$d$}\\
\hline
 & $0$ &  $0.1$ & $0.2$ & $0.3$ & $0.4$ & $0.5$ & $0.6$ & $0.7$ & $0.8$ & $0.9$ & $1.0$\\
\hline
$T_n$  &0.097 & 0.097 & 0.080 & 0.0896 &0.086 &0.072 &0.076 &0.077 &0.081 &0.075 &0.071\\
$\mbox{HvK}$  &0.146 & 0.138& 0.132 & 0.130 &0.124 &0.122 &0.110 &0.103 &0.102 & 0.099 &0.110\\
\mbox{Durot} & 0.021 & 0.019 & 0.018 & 0.013 & 0.015 & 0.018 &0.021 &0.024 &0.015 &0.018 &0.021\\
\mbox{Durot, $T_n$} & 0.003 & 0.003 & 0.004 & 0.001 & 0.002 & 0.003 &0.002 &0.001 &0.001 &0.001 &0.000\\
\hline
\end{tabular}
\end{table}

\begin{table}[ht!]
\centering
\caption{Estimated powers for model (\ref{F_d}), where $\a=0.1$, $n=50$, and $d=-1,-0.9,\dots,-0.1$. The estimation interval is $[0,F_0^{-1}(0.8)]$.}
\label{table:powers_d3}
\vspace{0.5cm}
\begin{tabular}{|l|c|c|c|c|c|c|c|c|c|c|c|}
\hline
 & \multicolumn{10}{c|}{$d$}\\
\hline
 & $-1$ &  $-0.9$ & $-0.8$ & $-0.7$ & $-0.6$ & $-0.5$ & $-0.4$ & $-0.3$ & $-0.2$ & $-0.1$\\
\hline
$T_n$  & 0.880 & 0.726 & 0.569 & 0.433 & 0.332 & 0.246 &0.204 &0.163 &0.140 &0.127\\
$\mbox{HvK}$ &0.965 & 0.864 & 0.766 & 0.686 &0.544 &0.483 &0.434 &0.326 &0.299& 0.279\\
\mbox{Durot} & 0.645& 0.524 & 0.399& 0.303 & 0.231 & 0.168 &0.127 &0.097 &0.080 &0.075\\
\mbox{Durot, $T_n$} & 0.742& 0.569 & 0.395& 0.253 & 0.181 & 0.121 &0.067 &0.063 &0.038 &0.030\\
\hline
\end{tabular}
\end{table}

\begin{table}[ht!]
\centering
\caption{Estimated rejection probabilities for model (\ref{F_d}) under the null hypothesis, where $\a=0.1$, $n=50$, and $d=0,0.1,\dots,1$. The estimation interval is $[0,F_0^{-1}(0.8)]$.}
\label{table:powers_d4}
\vspace{0.5cm}
\begin{tabular}{|l|c|c|c|c|c|c|c|c|c|c|c|c|}
\hline
 & \multicolumn{11}{c|}{$d$}\\
\hline
 & $0$ &  $0.1$ & $0.2$ & $0.3$ & $0.4$ & $0.5$ & $0.6$ & $0.7$ & $0.8$ & $0.9$ & $1.0$\\
\hline
$T_n$ & 0.101 & 0.103& 0.101 & 0.102& 0.096 & 0.094 &0.087 &0.091 &0.085 &0.073 &0.074\\
$\mbox{HvK}$  &0.256 & 0.229& 0.192 & 0.188 &0.170 &0.139 &0.145 &0.132 &0.121 & 0.131 &0.112\\
 \mbox{Durot} & 0.060 & 0.047 & 0.043 & 0.037 & 0.027 & 0.029 &0.037 &0.034 &0.028 &0.026 &0.025\\
 \mbox{Durot, $T_n$} & 0.024 & 0.019 & 0.016 & 0.009 & 0.013 & 0.009 &0.005 &0.006 &0.006 &0.004 &0.004\\
\hline
\end{tabular}
\end{table}

\begin{figure}[!ht]
\begin{center}
\includegraphics[scale=0.45]{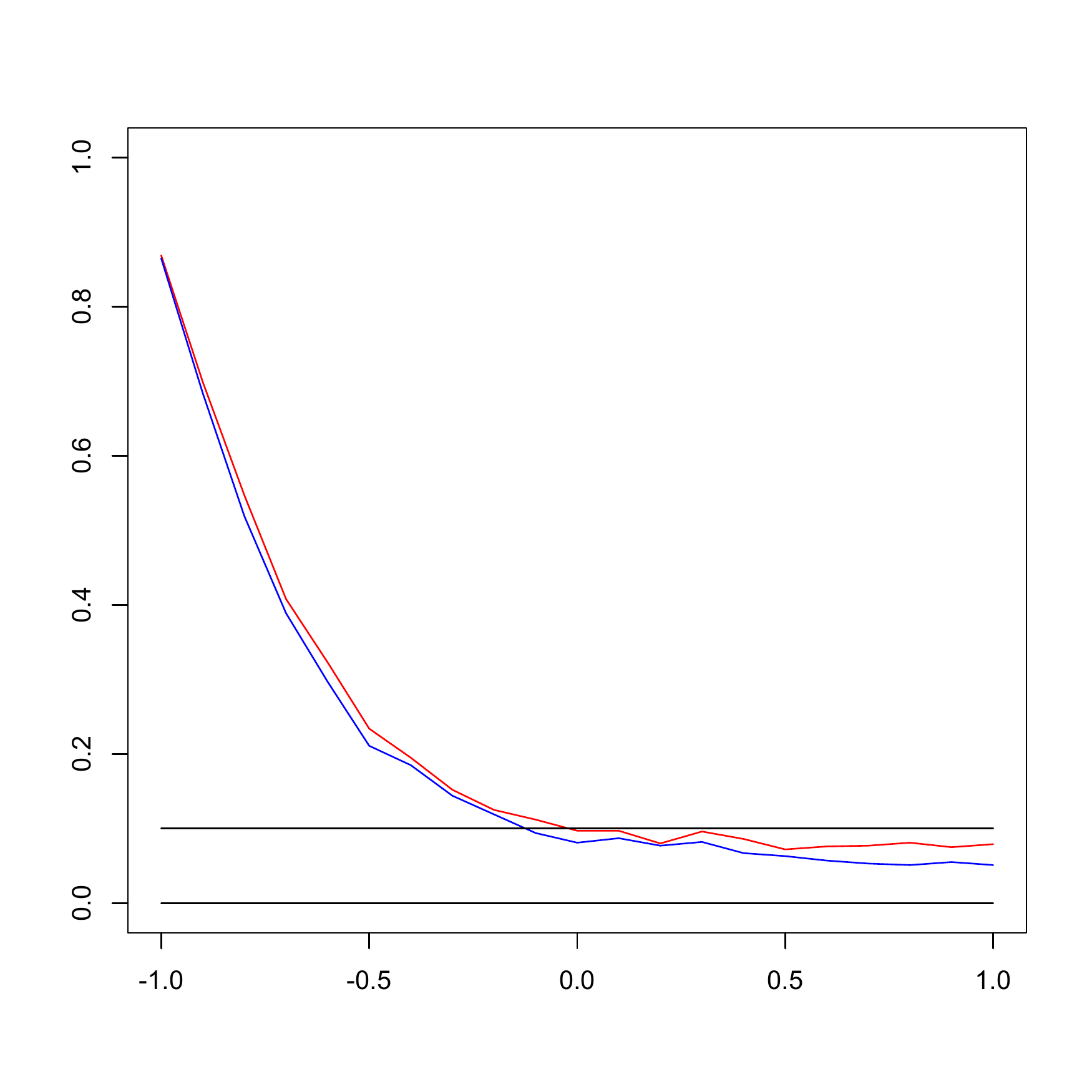}
\includegraphics[scale=0.45]{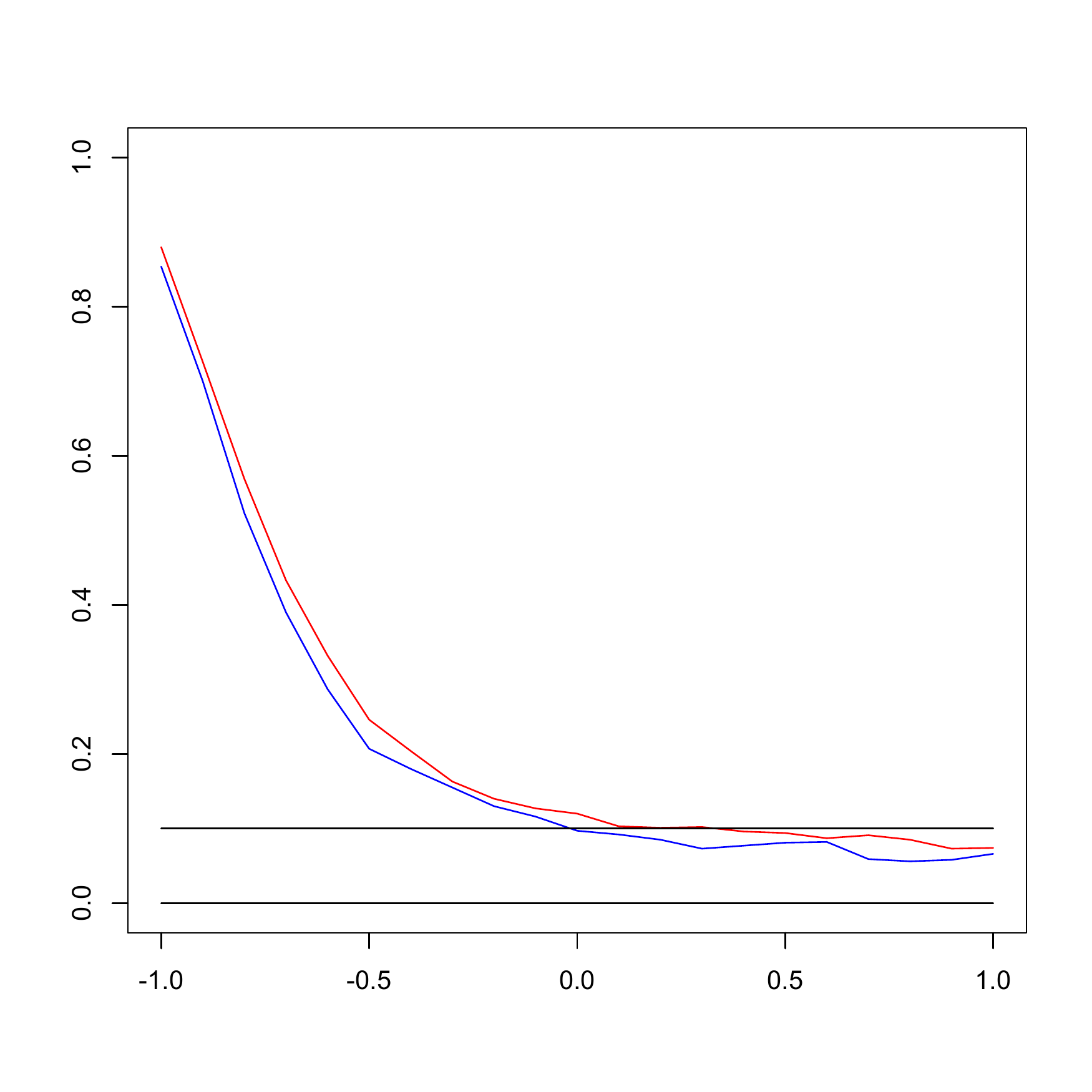}
\end{center}
\caption{The estimated power curves for the family $\{F^{(d)}:d\in[-1,1]\}$ and the isotonic test statistic $T_n$, defined by (\ref{eq:defTn}), for critical values estimated by bootstrapping from a penalized and smoothed isotonic estimate (red) and for critical values estimated by bootstrapping from the isotonic estimate itself (blue). The sample size $n=50$ and the estimation interval is $[0,F_0^{-1}(0.95)]$ in the left panel, $[0,F_0^{-1}(0.8)]$ in the right panel.}
\label{fig:powers50_naive}
\end{figure}

It is also of interest to compare the powers of the procedure, based on bootstrapping from a penalized and smoothed version of the hazard, with the powers obtained by just bootstrapping from the isotonic piecewise constant hazard estimate without any smoothing or penalizing. This is done in Figure \ref{fig:powers50_naive}, where it is seen that the difference is not very large for this family (and this sample size). The general trend is that bootstrapping from the isotonic estimate itself gives more conservative critical values.

\begin{figure}[!ht]
\begin{center}
\includegraphics[scale=0.65]{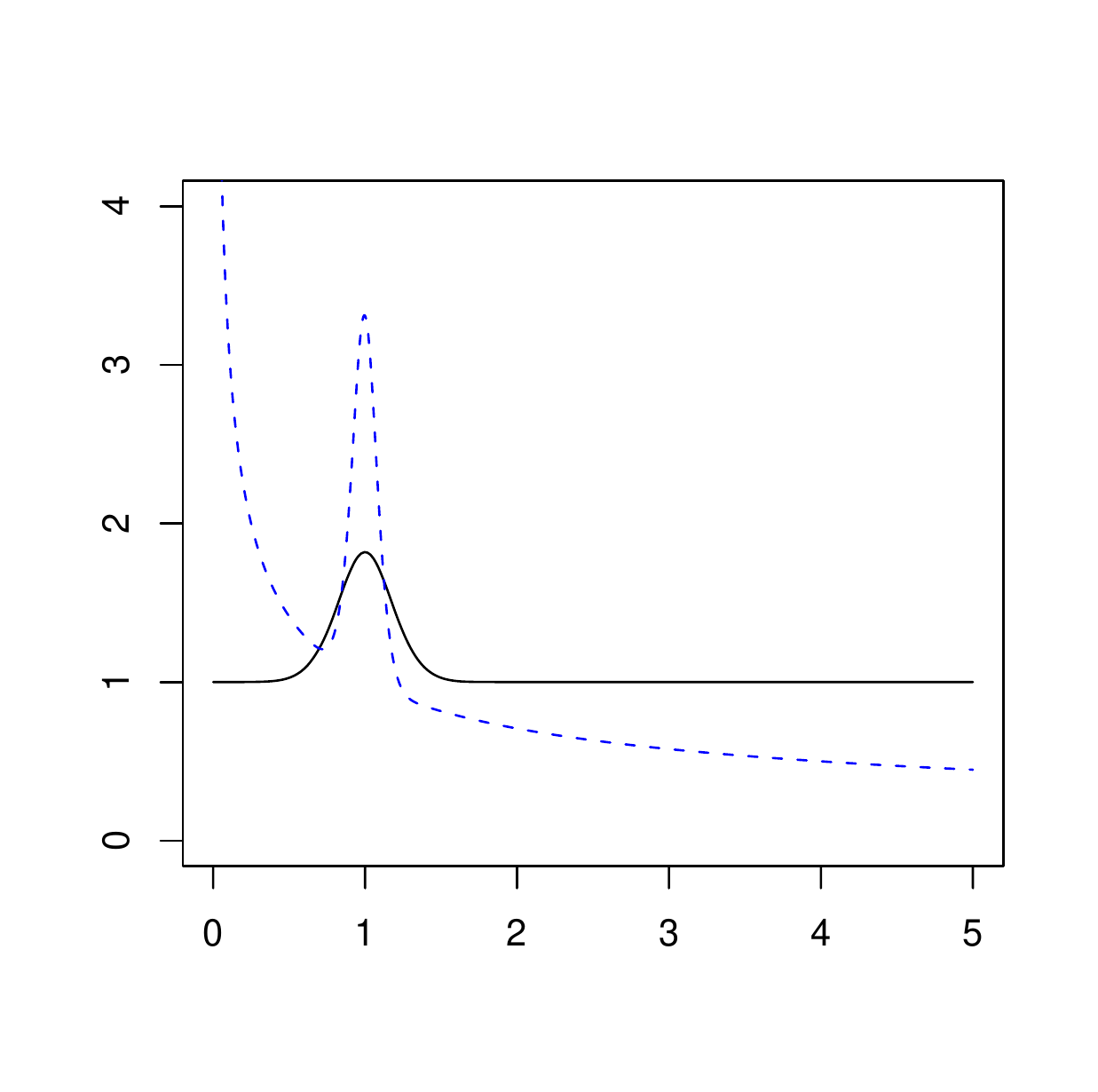}
\includegraphics[scale=0.65]{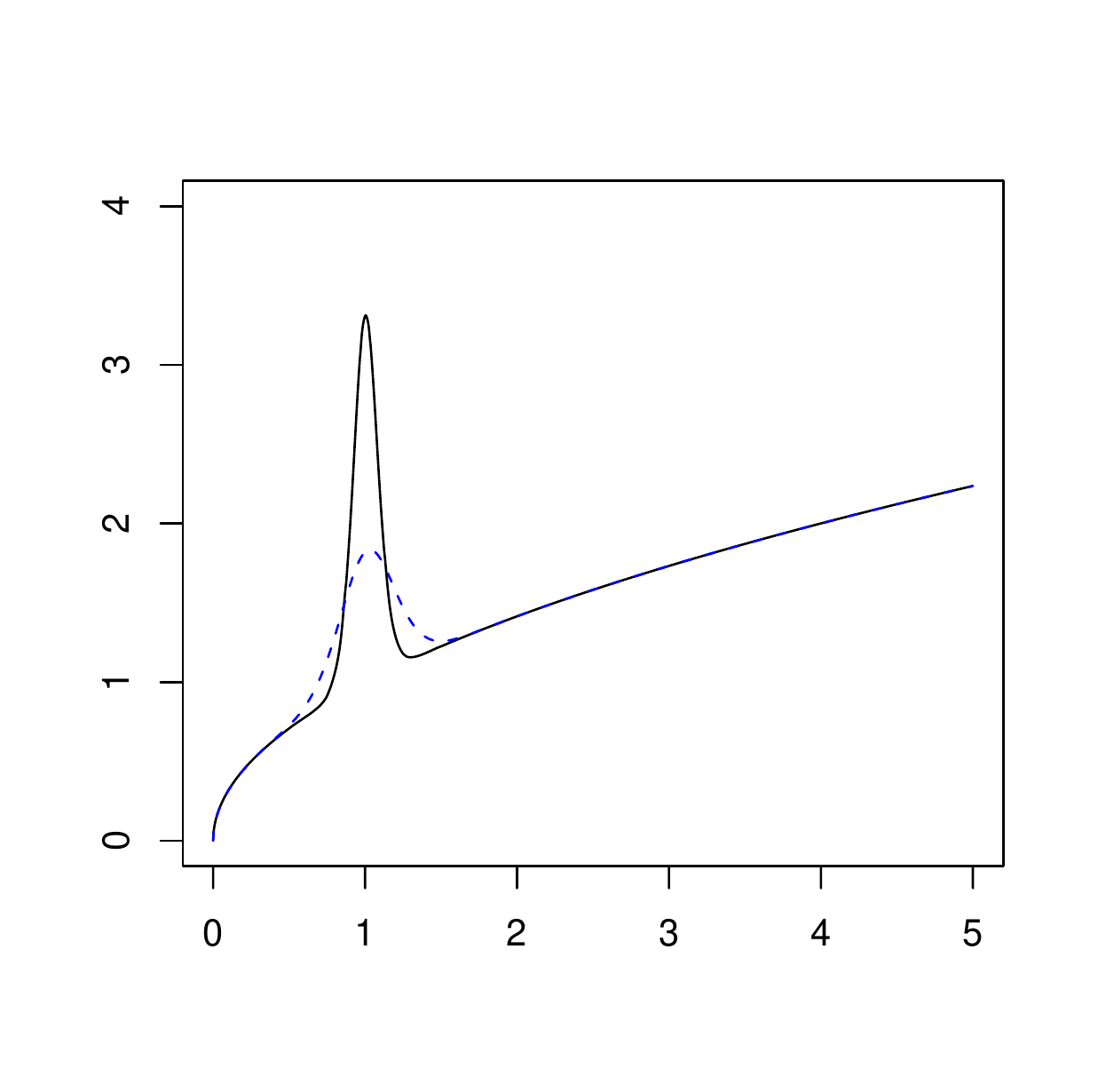}
\end{center}
\caption{The hazard function $h_{\b,\g,\m,\s}$. The solid line in the left panel corresponds to $\b=0.3,\,\g=0,\,\m=1$ and $\s=0.2$; the dashed line with $\b=0.3,\,\g=-0.5,\,\m=1$ and $\s=0.1$. In the right panel the solid line corresponds to $\b=0.3,\,\g=0.5,\,\m=1$ and $\s=0.1$; the dashed line with $h_{\b,\g,\m,\s}$, for $\b=0.3,\,\g=0.5,\,\m=1$ and $\s=0.2$.}
\label{fig:haz_gamma0}
\end{figure}

In \mycite{hallk:05} also the model, where the hazard function is of the form
\begin{equation}
\label{bump_family}
h_{\b,\g,\m,\s}(x)=x^{\g}\exp\left\{\b\left(2\pi\s^2\right)^{-1/2}\exp\left\{-(x-\m)^2/(2\s^2)\right\}\right\},
\end{equation}
is studied. Typical members of the family are shown in Figure \ref{fig:haz_gamma0}. For this family $T_n$ also seems to provide the most ``all-round" test, since it is more powerful than the HvK test for the global alternatives, where $\g<0$ (note that the hazard is globally decreasing for these alternatives) and more powerful for detecting the local disturbances where $\g>0$ than the \mycite{gijbheck:04} and \mycite{prospyke:67} tests. Note that the HvK test gives  rejection probabilities which are all above the $10\%$ level in the null hypothesis region. For the exponential distribution ($\b=\g=0$) the rejection probability is even close to $50\%$! The test, based on $T_n$ also gives a rejection probability which is too high here. This is probably caused by boundary effects and could possibly be remedied by adding heavier penalties in the cusum diagram at the beginning and end of the interval which is considered. The test was computed for the interval $[0,F_0^{-1}(0.95)]$.

The test based on $T_n$ also has higher power for this family than the Durot test, except when $\g=0$. When $\g=0$ (and $\b=0.3$) the hazard is constant except for a local bump (see Figure \ref{fig:haz_gamma0}), so the isotonic projection is the constant hazard. Since the critical values in the Durot  test are specifically based on a (locally) constant hazard, its behavior in this situation is not surprising, because resampling from the exponential distribution is in this case almost the same as resampling from the isotonic projection of the real hazard.

\begin{table}[!ht]
\centering
\caption{Estimated rejection probabilities for model (\ref{bump_family}), for $\a=0.1$ and $n=50$. The numbers in italics are rejection probabilities under the null hypothesis. The values for the HvK, PP and GH tests were taken from Table 1, p.\ 1124, in {\small\sc Hall and van Keilegom} (2005).}
\label{table:simstudy1}
\vspace{0.5cm}
\begin{tabular}{|l|c|c|c|c|c|c|}
\hline
 & \multicolumn{6}{c|}{$\gamma$}\\
\hline
 $\mbox{Parameter}$&  $\mbox{Test}$ & $-0.5$ &  $-0.25$ & $0$ & $0.5$ & $1$\\
\hline
$\beta=0$ & $T_n$ & 1.000 & 0.792 & {\it 0.213} &{\it 0.050} &{\it 0.076}\\
 & $\mbox{HvK}$ & 0.844 & 0.525 & {\it 0.437} & {\it 0.189}& {\it 0.121}\\
 & $\mbox{Durot}$ & 0.704 & 0.307 & {\it 0.096} & {\it 0.031}& {\it 0.028}\\
 & $\mbox{PP}$ & 1.00 & 0.800 & {\it 0.100} & {\it 0.000}& {\it 0.000}\\
 & $\mbox{GH}$ & 0.983 & 0.416 & {\it 0.100} & {\it 0.034}& {\it 0.027}\\
\hline
$\s=0.1$ & $T_n$ & 0.985 & 0.549 & 0.229& 0.497 & 0.585\\
$\m=1$ & $\mbox{HvK}$ & 0.675 & 0.753 & 0.772 & 0.656& 0.508\\
$\beta=0.3$ &  $\mbox{Durot}$ & 0.501 & 0.417 & 0.320 & 0.182& 0.107\\
 & $\mbox{PP}$ & 0.997 & 0.458 & 0.019 & 0.000& 0.000\\
 & $\mbox{GH}$ & 0.962 & 0.291 & 0.178 & 0.176& 0.154\\
\hline
$\s=0.2$ & $T_n$ & 0.991 &  0.605 & 0.172& 0.214 & 0.216\\
$\m=1$ & $\mbox{HvK}$ & 0.715 & 0.714 & 0.663& 0.443& 0.277\\
$\beta=0.3$ &  $\mbox{Durot}$ & 0.545 & 0.346& 0.218 & 0.090& 0.045\\
 & $\mbox{PP}$ & 0.999 & 0.588 & 0.053 & 0.000& 0.000\\
 & $\mbox{GH}$ & 0.968 & 0.301 & 0.114 & 0.065& 0.054\\
\hline
\end{tabular}
\end{table}

\section{Appendix}
\label{sec:appendix}
{\bf Proof of Lemma \ref{lem:consapp}:} Let $E_1,E_2,\ldots,E_n$ be an i.i.d.\ sequence of standard exponential random variables. Define
$$
X_i=H^{-1}(E_i),\,\,\, Y_i=H_{a,b}^{-1}(E_i)  \mbox{ for } 1\le i\le n.
$$
Then the $X_i$'s  and the $Y_i$'s are samples from the distributions with cumulative hazard $H$ and $H_{a,b}$ respectively.  Denote by $U_n$ the test statistic (\ref{eq:defTn}) based on the $Y_i$'s and by $V_n$ the statistic based on the $X_i$'s. Furthermore, define the  function $\phi\,:\,[a,b]\rightarrow [a,b]$ by
$
\phi(x)=H_{a,b}^{-1}(H(x)).
$
Note that $\phi$ is convex and increasing on $[a,b]$ and that $Y_i=\phi(X_i)\le X_i$ for all $i$. Moreover, using obvious notation,
$$
\F_n^X(x)=\frac1n\#\left\{i\,:\,X_i\le x\right\}=\frac1n\#\left\{i\,:\,\phi(X_i)\le \phi(x)\right\}=\frac1n\#\left\{i\,:\,Y_i\le \phi(x)\right\}=\F_n^Y(\phi(x)).
$$
Consequently, also $\H_n^X(x)=\H_n^Y(\phi(x))$, where these functions refer to the empirical cumulative hazards based on the sample of $X_i$'s and $Y_i$'s respectively.
Now define $\bar H_n(x)=\hat H_n^Y(\phi(x))$, where the latter denotes the greatest convex minorant of the  empirical hazard function based on the $Y_i$'s, evaluated at $\phi(x)$.
Then $\bar H_n$ is  {\it a minorant} of $\H_n^X$, i.e.,
$$
\bar H_n(x)=\hat H_n^Y(\phi(x))\le \H_n^Y(\phi(x))=\H_n^X(x).
$$
Moreover, it is also convex. Indeed, using monotonicity and convexity of $\hat H_n^Y$ and convexity of $\phi$ we have for $\alpha\in(0,1)$ and $x,y\in[a,b]$
\begin{eqnarray*}
\bar H_n(\alpha x+(1-\alpha)y)&=&\hat H_n^Y(\phi(\alpha x+(1-\alpha)y))\le\hat H_n^Y(\alpha\phi(x)+(1-\alpha)\phi(y))\\
&\le&\alpha \hat H_n^Y(\phi(x))+(1-\alpha)\hat H_n^Y(\phi(y))=\alpha \bar H_n(x)+(1-\alpha)\bar H_n(y).
\end{eqnarray*}
Hence, the convex minorant $\bar H_n$ of $\H_n^X$ is smaller than or equal to the greatest convex minorant $\hat H_n^X$ of $\H_n^X$:
$$
\bar H_n(x)\le\hat H_n^X(x)\le \H_n^X(x)\Rightarrow  \F_n^X(x)-\hat F_n^X(x)\le  \F_n^X(x)-\bar F_n(x),
$$
where we use the obvious notation relating cumulative hazards to distribution functions. This implies that
\begin{eqnarray*}
U_n:&=&\int_{[a,b]}\left(\F_n^Y(x-)-\hat F_n^Y(x)\right)d\F_n^Y(x)=\int_{[a,b]}\left(\F_n^Y(\phi(x)-)-\hat F_n^Y(\phi(x))\right)d\F_n^Y(\phi(x))\\
&=&\int_{[a,b]}\left(\F_n^X(x-)-\bar F_n(x)\right)d\F_n^X(x)\ge \int_{[a,b]}\left(\F_n^X(x-)-\hat F_n^X(x)\right)d\F_n^X(x)=:V_n.
\end{eqnarray*}
Noting that $P_{H}\left(T_n\ge t\right)=P(V_n\ge t)$  and $P_{H_{a,b}}\left(T_n\ge t\right)=P(U_n\ge t)$, the result follows.
\hfill$\Box$

\vspace{0.5cm}
\noindent
{\bf Proof of Theorem \ref{th:bootstrap1}}. The result follows from Theorem \ref{th:limitdist}, if we can show that the estimate $\tilde H_n$, which generates the bootstrap samples, has the property that the corresponding estimates $\tilde f_n$, $\tilde h_n$ and $\tilde h_n'$ of $f_0$, $h_0$ and $h_0'$, respectively, will be consistent (in an almost sure sense), since in this case the integrals
\begin{equation}
\label{mu_bootstrap}
\int_a^b \,\left(\frac{2\tilde h_n(t)\tilde f_n(t)}{\tilde h_n'(t)}\right)^{1/3} \tilde f_n(t)\,dt,
\end{equation}
and
\begin{equation}
\label{sigma_bootstrap}
\int_a^{b}\left(\frac{2\tilde h_n(t)\tilde f_n(t)}{\tilde h_n'(t)}\right)^{4/3} \tilde f_n(t)\,dt
\end{equation}
will converge (almost surely) to the integrals defining (\ref{mu}) and (\ref{sigma}). Moreover, the consistency will ensure that the derivative $\tilde h_n'$ will stay away from zero for large $n$, implying that the distribution, generating the bootstrap samples will satisfy the conditions of Theorem \ref{th:limitdist} for large $n$, implying that the asymptotic normality result also holds for the test statistics, computed for the bootstrap samples (with parameters $\mu^*_n$ and $\sigma^*_n$, derived from (\ref{mu_bootstrap}) and (\ref{sigma_bootstrap})).

But the uniform consistency of the estimates $\tilde h_n$ and $\tilde f_n$ on $[a,b]$ is proved in \mycite{geurt_piet:11b} (here we also use the penalization of $\hat h_n$!), and the consistency of $\tilde h_n'$ on the interior of $[a,b]$ is ensured by the choice of the bandwidth $cn^{-1/4}$. Since this choice of bandwidth also ensures that the right limit of $\tilde h_n'$ at $a$ and the left limit of $\tilde h_n'$ at $b$ will be positive for all large $n$, we indeed have:
$$
\int_a^b\left(\frac{2\tilde h_n(t)\tilde f_n(t)}{\tilde h_n'(t)}\right)^{1/3} \tilde f_n(t)\,dt
\stackrel{a.s.}\longrightarrow \int_a^b f_0(t)\,\left(\frac{2h_0(t)f_0(t)}{h_0'(t)}\right)^{1/3} dF_0(t)
$$
and
$$
\int_a^{b}\left(\frac{2\tilde h_n(t)\tilde f_n(t)}{\tilde h_n'(t)}\right)^{4/3} \tilde f_n(t)\,dt
\stackrel{a.s.}\longrightarrow \int_a^b f_0(t)^2\,\left(\frac{2h_0(t)f_0(t)}{h_0'(t)}\right)^{4/3} dF_0(t),
$$
and the result now follows.\eop

\vspace{0.3cm}
\begin{remark}
{\rm In order to get a good estimate of the critical value, one wants to choose a small bandwidth  in estimating the hazard function for the bootstrap samples, in order to minimize the bias. However, since one also wants to estimate the derivative $h_0'$ consistently on the interval $[a,b]$, the bandwidth cannot be too small. As an example, the choice $b_n=n^{-1/3}$ is too small for this purpose.
This motivated the choice of the bandwidth of order $n^{-1/4}$, but also other choices are possible.
\mycite{silverman:78} gives the necessary and sufficient condition:
$$
\frac{n b_n^3}{\log(1/b_n)}\to\infty,
$$
for uniform consistency of a kernel estimate of the density in ordinary density estimation (see his Theorem C, p.\ 182), where $b_n$ again denotes the bandwidth.
}
\end{remark}

\vspace{0.3cm}
\noindent
{\bf Acknowledgement}
The first author wants to thank David Steinsaltz for introducing him to the topic of this paper.

\end{document}